\documentclass[10pt,twoside]{amsart}
\usepackage{lmodern}
\usepackage[T1]{fontenc}
\usepackage{amsmath,amsfonts,amsthm,amssymb,amscd}
\usepackage{mathrsfs}
\usepackage{latexsym}
\usepackage[pdftex]{graphicx}
\usepackage{enumitem}
\usepackage{bbm}
\usepackage{accents}

\usepackage{geometry}
\newtheorem{lemma}{Lemma}[section]
\newtheorem{theorem}{Theorem}[section]

\usepackage{hyperref}  
\hypersetup{colorlinks=true, linktoc=all,
    linkcolor=blue,
		citecolor=blue,
		urlcolor=cyan,
		bookmarksopen=true
		}
\usepackage{hyperref}  

\usepackage[english]{babel}

\title{Distribution of products of shifted primes in arithmetic progressions with increasing difference}
\author{Zarullo Rakhmonov}
\address{A.Dzhuraev Institute of Mathematics,  National Academy of Sciences of Tajikistan}
\email{zarullo.rakhmonov@gmail.com, zarullo-r@rambler.ru}
\date{}

\begin{document}

\begin{abstract}
We obtain an asymptotic formula for the number of primes $p\leq x_1$, $p\leq x_2$ such that $p_1(p_2+a)\equiv l \pmod q$ with $q\leq x^{\ae_0}$, $x_1\geq x^{1-\alpha}$, $x_2\geq x^{\alpha}$,
$$
{\kappa}_0=\frac{1}{2.5+\theta+\varepsilon}, \quad \alpha\in \left[(\theta+\varepsilon)\frac{\ln q}{\ln x}, 1-2.5\frac{\ln q}{\ln x}\right],
$$
where $\theta=1/2$, if $q$ is a cube free and $\theta=\frac{5}{6}$ otherwise. This is the refinement and generalization of the well-known formula of A.~A.~Karatsuba.
Bibliography: 39 references.
\end{abstract}

\maketitle

\section{Introduction}
For the Dirichlet character $\chi$ modulo $q$, the Chebyshev function is defined by the equality
$$
\psi(y,\chi)=\sum_{n\le y}\Lambda(n)\chi(n),
$$
where $\Lambda(n)$ is the von Mangoldt function. We will also use the following notation: $\varphi(q)$ is Euler`s totient function, $\mu(n)$ is the Mobius function, and $\mathscr{L}=\ln q$.

It is known that, under the assumption of the Extended Riemann Hypothesis, the following estimates hold:
\begin{align}
&t(x;q)=\sum_{\chi\hspace{-9pt}\mod\hspace{-2pt}q}\max_{y\leq x}|\psi(y,\chi)|\ll x+x^\frac12q(\ln xq)^2,\label{formula-otsenka t(x;q) v pred RGR}\\
&T(x;Q)=\sum_{q\le Q}\frac{q}{\varphi(q)}\sideset{}{^*}\sum_{\chi}\max_{y\leq x}|\psi(y,\chi)|\ll x^\frac12Q^2(\ln xQ)^2, \label{formula-otsenka T(x;Q) v pred RGR}
\end{align}
where the symbol $^*$ indicates that the summation is taken over all primitive characters modulo $q$.

In solving a number of problems in the theory of prime numbers, it is sufficient that the estimates for $T(x;Q)$ and $t(x;q)$ are close to those given in (\ref{formula-otsenka t(x;q) v pred RGR}) and (\ref{formula-otsenka T(x;Q) v pred RGR}).

A.A.~Karatsuba~\cite{Karatsuba-DANSSSR-1970-192-4,KaratsubaAA-2008} developed a method for solving ternary multiplicative problems, by means of which he obtained an estimate for the simplest case of the quantity $t(x;q)$ and solved the problem of the distribution of numbers of the form $p_1(p_2 + a)$ in arithmetic progressions with growing difference in the following formulation.

{\theorem\label{Teorema-AAK-RPSdPrChVArProgr}
Let $\varepsilon \in \left(0, \frac{1}{4}\right]$; $x \ge x_0(\varepsilon)$ be a sufficiently large positive number; let $q$ be a prime number such that $q \le x^{\ae_0}$, where $\kappa_0 = 1/(4.6 + \varepsilon)$; $(a, q) = 1$, $(l, q) = 1$; and let $\alpha$ be an arbitrary number in the interval
$$
\left[\left(\frac{1}{2} + \varepsilon\right)\frac{\mathscr{L}}{\ln x}, 1 - 4.1\frac{\mathscr{L}}{\ln x}\right],
$$
where $x_1 \ge x^{1 - \alpha}$, $x_2 \ge x^{\alpha}$; $p_1$, $p_2$ are prime numbers; and let $\pi_2(x_1,x_2,a,l)$ denote the number of pairs with $p_1 \le x_1$, $p_2 \le x_2$ such that $p_1(p_2 + a) \equiv l \hspace{-2pt} \pmod q$; and let $\delta > 0$ be an arbitrarily small number. Then the following asymptotic formula holds:
$$
\pi_2(x_1,x_2,a,l) = \frac{\pi(x_1)\pi(x_2)}{\varphi(q)} + O\left(\left(x_1x_2\right)^{1+\delta}q^{-1 - \frac{\varepsilon^2}{1024}}\right),
$$
where the constant in the $O$-term depends only on $\varepsilon$.}

A.A.~Karatsuba~\cite{Karatsuba-DANSSSR-1970-192-4}, in this work, noted that:
\begin{itemize}
\item The problem of the distribution in arithmetic progressions of numbers of the form $(p_1^n + a)f(p_2)$, where $p_1$ and $p_2$ are prime numbers and $f$ is a polynomial with integer coefficients, can be studied in essentially the same way;
\item The same method can also be applied to problems concerning the distribution of prime numbers in arithmetic progressions ``on average'' and to other related problems.
\end{itemize}

M.~M.~Petechuk~\cite{PetechukMM-IzvFNSSSR-1979-43-4}, applying A.~A.~Karatsuba's method of solving ternary multiplicative problems and using estimates for short character sums obtained by V.N.~Chubarikov~\cite{Chubarikov-Vestnik-MGU-1973-2}, proved an asymptotic formula for the sum
$$
S=\sum_{\substack{n\le x\\ n\equiv l\bmod q}}\tau_k(n),
$$
where $q = p^m$, $p$ is a fixed prime, $(l, q) = 1$, and $q \le x^{\frac{3}{8} + \varepsilon}$. Subsequently, A.~A.~Karatsuba and M.~M.~Petechuk obtained an asymptotic formula for the sum $S$ in the case where $q$ is a prime number and $q \le x^{\frac{4}{k} - \varepsilon}$, with $k \ge k_0 \ge 7$ (1979, presented at the seminar on analytic number theory at Moscow State University). The use of estimates for short character sums modulo cube-free moduli allowed H.~Iwaniec and J.~B.~Friendlander~\cite{Friendlander+Iwaniec-ActaArith-1985-45-3} to extend this result to the case of cube-free moduli.

In 1989, the author~\cite{RakhmonovZKh-Izv-ANSSSR-1989}, relying on A.A.~Karatsuba's method, gave an elementary proof of the estimate
\[
t(x;q) \ll (x + x^{\frac56}q^{\frac12} + x^{\frac12}q)x^\varepsilon.
\]
Using the same method, Pan Chen Dong and Pan Chen Biao~\cite{PanChD+HfnChB-OATCh-1991-Pekin} showed that
\[
T(x;Q) \ll (x + x^{\frac56}Q + x^{\frac12}Q^2)(\ln xQ)^4.
\]
As a consequence of this estimate, one obtains a theorem on the distribution of prime numbers in arithmetic progressions ``on average,'' the possibility of which was indicated in~\cite{Karatsuba-DANSSSR-1970-192-4}.

G.~Montgomery~\cite{Montgomeri-1974}, using his density theorem on the zeros of Dirichlet $L$-functions (the proof of which relies on the large sieve method), established the following estimates:
\begin{align*}
&t(x;q) \ll (x + x^{\frac57}q^{\frac57} + x^{\frac12}q)(\ln xq)^{16},\\
&T(x;Q) \ll (xQ^{\frac23} + x^{\frac12}Q^2)(\ln xQ)^{11}.
\end{align*}
This result was refined by R.~Vaughan~\cite{Vaughan-1975}, who, using the large sieve method and a special representation for the logarithmic derivative of $L$-functions, proved that
\begin{align}
&t(x;q) \ll x(\ln xq)^3 + x^{\frac34}q^{\frac58}(\ln xq)^{\frac{23}8} + x^{\frac12}q(\ln xq)^{\frac72}, \label{formula-otsenka t(x;q)-Vaughan}\\
&T(x;Q) \ll x(\ln xQ)^3 + x^{\frac34}Q^\frac54(\ln xQ)^{\frac{23}8} + x^{\frac12}Q^2(\ln xQ)^{\frac72}. \label{formula-otsenka T(x;Q)-Vaughan}
\end{align}
The author~~\cite{RakhmonovZKh-IzvRAN-1993-57-4,RakhmonovZKh-DAN-1993-331-3,RZKh+NOO-ChebSb-2021-22-59(81),RakhmonovZKh-1994-IzvRAN, RakhmonovZKh-DAN-1996},  applying A.~A.~Karatsuba's method for solving ternary multiplicative problems in combination with a new analytic version of I.~M.~Vinogradov's method for estimating trigonometric sums with prime numbers, and using the approach of N.~M.~Timofeev~\cite{Timofeev-Izv-ANSSSR-1987}, where he studies the distribution of arithmetic functions in short intervals on average over progressions, followed by the application of Montgomery's theorem~\cite{Montgomeri-1974} on the fourth moment of Dirichlet $L$-functions, proved that
\begin{align}
&t(x;q) \ll x(\ln xq)^3 + x^{\frac45}q^{\frac12}(\ln xq)^{34} + x^{\frac12}q(\ln xq)^{34}, \label{formula-otsenka t(x;q)-RZKh}\\
&T(x;Q) \ll x(\ln xQ)^3 + x^{\frac45}Q(\ln xQ)^{34} + x^{\frac12}Q^2(\ln xQ)^c, \label{formula-otsenka T(x;Q)-RZKh}
\end{align}
where
\[
c = \begin{cases}
34, & \text{if } Q \le x^\frac53(\ln x)^{-\frac53}, \\
\frac{7}{2}, & \text{if } Q > x^{\frac56}(\ln x)^{-\frac56}.
\end{cases}
\]

Note that the estimates~\eqref{formula-otsenka t(x;q)-RZKh} and~\eqref{formula-otsenka T(x;Q)-RZKh} are more precise than~\eqref{formula-otsenka t(x;q)-Vaughan} and~\eqref{formula-otsenka T(x;Q)-Vaughan}, respectively, in the range
\begin{align*}
&x^\frac25(\ln x)^\frac15<q\le x^{\frac23}(\ln x)^{-\frac53}, \quad x^\frac15<Q\le x^\frac13,
\end{align*}
while for other values of $q$ and $Q$, the estimates coincide up to powers of $\ln xq$ and $\ln xQ$.

We now restate the estimate~\eqref{formula-otsenka t(x;q)-RZKh} in a form that will be used below in the proof of Theorem~\ref{Teorema-RZKh-RPSdPrChVArProgr}. We have:

{\lemma\label{Sled-teorema-RZKh-otsenka-t(x;q)}
Let $x \ge 2$ and $q \ge 1$. Then the following estimate holds:
\[
\sum_{\chi \bmod q} \max_{y \le x} \left| \sum_{p \le y} \chi(p) \right| \ll x\mathscr{L}_q^2 + x^{\frac45}q^{\frac12} \mathscr{L}_q^{33} + x^{\frac12}q\mathscr{L}_q^{33}.
\]
}

From estimate~\eqref{formula-otsenka T(x;Q)-RZKh}, the Bombieri-Vinogradov theorem on the ``average'' distribution of prime numbers in arithmetic progressions follows in the following form, which will also be used in the proof of Theorem~\ref{Teorema-RZKh-RPSdPrChVArProgr}:

{\lemma\label{Sledst-teorema-RZKh-otsenka-T(x;Q)}
Let $A$ be an arbitrary positive number. Then the following estimate holds:
\[
\sum_{q \le \sqrt{x}(\ln x)^{-A - 3.5}} \max_{y \le x} \max_{(l, q) = 1} \left| \pi(x; q, l) - \frac{\mathrm{Li}(x)}{\varphi(q)} \right| \ll \frac{x}{(\ln x)^{A+1}}.
\]
}

A.~A.~Karatsuba solved the problem of the distribution of numbers of the form $p_1(p_2 + a)$ in arithmetic progressions with growing difference, relying essentially on his estimate of the sum of a non-principal character modulo $q$ over a sequence of shifted prime numbers.

The problem of the distribution of values of a non-principal character on sequences of shifted prime numbers was first studied by I.~M.~Vinogradov~\cite{VinigradovIM-1938}. Using his method of estimating trigonometric sums over prime numbers, in 1938 he proved the following: \emph{if $q$ is an odd prime, $(a,q)=1$, and $\chi(a)$ is a non-principal character modulo $q$, then}
\[
T(\chi,x) = \sum_{p \le x} \chi(p + a) \ll x^{1 + \varepsilon} \sqrt{\frac{1}{q} + \frac{q}{\sqrt[3]{x}}}.
\]
In 1943, I.~M.~Vinogradov~\cite{VinigradovIM-1943} refined this estimate by proving that
\begin{equation}\label{Formula 0tsenka IMV-1943}
|T(\chi,x)| \ll x^{1 + \varepsilon} \left( \sqrt{\frac{1}{q} + \frac{q}{x}} + x^{-\frac16} \right).
\end{equation}
For $x \gg q^{1 + \varepsilon}$, this estimate is nontrivial and implies an \emph{asymptotic formula for the number of quadratic residues (non-residues) $\bmod\ q$ of the form $p + a$, with $p \le x$}. M.~Jutila~\cite{Jutila}, using estimate~\eqref{Formula 0tsenka IMV-1943}, showed that if $q$ is an odd prime, then
\[
G(q, l) \ll q^{\frac{11}8+\varepsilon},
\]
where $G(q,l)$ denotes the least Goldbach number in the arithmetic progression with difference $q$ and initial term $l$. A \emph{Goldbach number} is a number representable as the sum of two odd prime numbers.

Subsequently, I.~M.~Vinogradov obtained a nontrivial estimate for $T_1(\chi, x)$ when $x \ge q^{0.75 + \varepsilon}$, where $q$ is a prime number~\cite{VinigradovIM-1952,VinigradovIM-1953}. This result was unexpected. The point is that $T_1(\chi, x)$ can be expressed as a sum over the zeros of the corresponding Dirichlet $L$-function. Then, assuming the validity of the Generalized Riemann Hypothesis (GRH), one obtains a nontrivial estimate for $T_1(\chi, x)$ only when $x \ge q^{1+\varepsilon}$.

It seemed that this result contradicted the expected bounds. Yu.~V.~Linnik~\cite{Linnik-1973} wrote on this topic in 1971: \emph{``The investigations of I.~M.~Vinogradov in the asymptotics of Dirichlet characters are very important. As early as 1952, an estimate for the sum of Dirichlet characters over shifted primes $T_1(\chi)$ was obtained, which gave a power-saving in $x$ already when $x > q^{0.75 + \varepsilon}$. This estimate is of fundamental significance, since it surpasses, in depth, what is given by a direct application of the extended Riemann Hypothesis, and, apparently, in this direction represents a deeper truth than the hypothesis itself (assuming the hypothesis is true). Recently, this estimate was improved by A.~A.~Karatsuba.''}

In 1968, A.~A.~Karatsuba~\cite{KaratsubaAA-2008, KaratsubaAA-1968, KaratsubaAA-1970-1} developed a method that enabled him to obtain nontrivial estimates for short character sums in finite fields of fixed degree. In 1970, by advancing this method and combining it with I.~M.~Vinogradov's approach, he~\cite{KaratsubaAA-2008, KaratsubaAA-1970-2, KaratsubaAA-1970-3} proved the following result: \emph{if $q$ is a prime, $\chi(a)$ is a non-principal character modulo $q$, and $x \ge q^{1/2 + \varepsilon}$, then}
\begin{equation}\label{formula-otseka-AAK-T(chi)}
T(\chi, x) \ll xq^{-\frac{\varepsilon^2}{1024}}.
\end{equation}
A.~A.~Karatsuba applied these estimates to derive asymptotic formulas for the number of quadratic residues and non-residues of the form $p + a$, and for the number of products of shifted primes of the form $p_1(p_2 + a)$ in arithmetic progressions with increasing difference~\cite{Karatsuba-DANSSSR-1970-192-4, KaratsubaAA-2008}.

The author generalized the estimate~(\ref{Formula 0tsenka IMV-1943}) to the case of a composite modulus and proved the following statement~\cite{RakhmonovZKh-1986-UMN, RakhmonovZKh-1986-DANRT, RakhmonovZKh-1994-TrMIRAN}: \emph{Let $D$ be a sufficiently large natural number, $\chi$ a non-principal character modulo $D$, $\chi_q$ the primitive character inducing $\chi$, and $q_1$ the product of the prime divisors of $D$ that do not divide $q$. Then}
\begin{equation}\label{formula-RZKh-Otsenka-T(chi)-dlinn-summa}
T(\chi, x) \ll x \ln^4 x \left(\sqrt{\frac{1}{q} + \frac{q}{x} \tau^2(q_1)} + x^{-\frac{1}{6}} \tau(q_1)\right) \tau(q).
\end{equation}
Using this estimate, he also proved~\cite{RakhmonovZKh-1986-UMN, RakhmonovZKh-1986-IzvANRT} that for a sufficiently large odd natural number $D$, the following estimate holds:
\begin{equation}\label{formula-RZKh-otsen-naim-gold-chisla}
G(D, I) \ll D^{c + \varepsilon},
\end{equation}
where $\varepsilon$ is an arbitrarily small positive constant, and $c$ is the infimum of those numbers $a$ for which, for some constant $A > 2$,
$$
\sum_{\chi \bmod D} N(\alpha, T, \chi) \ll (DT)^{2a(1 - \alpha)} (\ln DT)^A.
$$
From Huxley's ``density'' theorem~\cite{Huxley-1942-Inv.Math}, it follows that for $A = 14$ in the formula above, one has $c \le \frac{6}{5}$.

In 2010, J.B.~Friedlander, K.~Gong, and I.E.~Shparlinski~\cite{Friedlander-Gong-Shparlinski} showed that for a composite modulus $q$, a nontrivial estimate for the sum $T(\chi_q, x)$ exists even when $x$---the length of the sum---is of order smaller than $q$. They proved the following: \emph{For any primitive character $\chi_q$ and every $\varepsilon > 0$, there exists $\delta > 0$ such that for all $x \ge q^{\frac{8}{9} + \varepsilon}$, the following estimate holds:}
\begin{equation}\label{Formula-otsenka-Friedlander-Gong-Shparlinski}
T(\chi_q, x) \ll x q^{-\delta}.
\end{equation}

In 2013, the author~\cite{RakhmonovZKh-2013-DANRT, RakhmonovZKh-2013-IzvSarUniv, RakhmonovZKh-2014-Ch.sbor-15-2(50)} proved that \emph{if $q$ is a sufficiently large natural number, $\chi_q$ is a primitive character modulo $q$, and $\varepsilon > 0$ is an arbitrarily small positive constant, then for $x \ge q^{\frac{5}{6} + \varepsilon}$ we have}
\begin{align*}
T(\chi_q, x) \ll x \exp\left(-\sqrt{\mathscr{L}}\right).
\end{align*}
In 2021, Bryce Kerr~\cite{Kerr-2021} proved the estimate~(\ref{Formula-otsenka-Friedlander-Gong-Shparlinski}) for $x \ge q^{\frac{3}{4} + o(1)}$.

As already mentioned above, nontrivial estimates of the sum $T(\chi, x)$, where $\chi$ is a non-principal character modulo a prime $q$, have been applied in problems on the distribution of products of shifted primes and on the least Goldbach numbers in short arithmetic progressions. When solving these problems for a composite modulus $q$, alongside nontrivial estimates for $T(\chi, x)$ with primitive characters, analogous estimates are also required for derived characters.

Therefore, it is natural to consider the problem of obtaining nontrivial estimates for sums of the form $T(\chi, x)$, where $\chi$ is a non-principal character modulo a composite number $q$, not only for long sums (estimate~(\ref{formula-RZKh-Otsenka-T(chi)-dlinn-summa})), that is, when the length of the sum is of order greater than the modulus of the character $\chi$, but also for short sums, when $x$---the length of the sum---is of order smaller than $q$. The author obtained the following nontrivial estimates.

{\lemma\label{Teorema-RZKh-Otsenka-T(chi)-q-bezcub} {\rm \cite{RakhmonovZKh-2017-TrMIRAN-299, RakhmonovZKh-2017-DANRT-9}.}
Let $q$ be a sufficiently large natural number, $\chi$ a non-principal character modulo $q$, and let $\chi_d$ be the primitive character modulo $d$ induced by $\chi$, where $d$ is cube-free. Suppose $(a, q) = 1$ and $\varepsilon > 0$ is an arbitrarily small constant. Then for $x \ge q^{\frac{1}{2} + \varepsilon}$, we have
\begin{align*}
T(\chi, x) \ll x \exp\left(-0.6 \sqrt{\mathscr{L}}\right),
\end{align*}
where the implied constant depends only on $\varepsilon$.}

{\lemma\label{Teorema-RZKh-Otsenka-T(chi)-q-sost} {\rm \cite{RakhmonovZKh-2025-SummNeglChar-arxiv, RakhmonovZKh-2018-Springer}.}
Let $q$ be a sufficiently large natural number, $\chi$ a non-principal character modulo $q$, and $\chi_d$ the primitive character modulo $d$ induced by $\chi$. Assume $(a, q) = 1$, $\varepsilon > 0$ is arbitrarily small, and $q > \exp \sqrt{2 \ln D}$. Then for $x \ge q^{\frac{5}{6} + \varepsilon}$, we have
\begin{align*}
T(\chi, x) \ll x \exp\left(-0.6 \sqrt{\mathscr{L}}\right),
\end{align*}
where the implied constant depends only on $\varepsilon$.}

A.~A.~Karatsuba in his paper~\cite{Karatsuba-DANSSSR-1970-192-4} also noted that in Theorem~\ref{Teorema-AAK-RPSdPrChVArProgr}, the upper bound for $q$ in the case where $q$ is a prime number can be significantly increased, but not beyond $x^{\ae_1}$. That is, the value
\[
\kappa_0 = \frac{1}{4.6 + \varepsilon}
\]
can be replaced by
\[
\kappa_1 = \frac{1}{2.5 + \varepsilon},
\]
which follows from the conditional estimate~(\ref{formula-otsenka t(x;q) v pred RGR}).

In the present work, the author was able to prove Karatsuba's theorem on the distribution of numbers of the form $p_1(p_2 + a)$ in arithmetic progressions with growing difference, by using a new estimate for the average values of the Chebyshev functions over all Dirichlet characters of a given modulus (Lemma~\ref{Sled-teorema-RZKh-otsenka-t(x;q)}), and nontrivial estimates for short sums of values of non-principal characters modulo $q$ in the sequence of shifted primes (Lemmas~\ref{Teorema-RZKh-Otsenka-T(chi)-q-bezcub} and~\ref{Teorema-RZKh-Otsenka-T(chi)-q-sost}), under the condition
\begin{equation}\label{formula-opred-ae0}
\kappa_0 = \frac{1}{2.5 + \theta + \varepsilon}, \qquad \theta = \left\{
  \begin{array}{ll}
    \dfrac{1}{2}, & \text{if $q$ is cube-free;} \vspace{8pt} \\
    \dfrac{5}{6}, & \text{otherwise.}
  \end{array}
\right.
\end{equation}

{\theorem\label{Teorema-RZKh-RPSdPrChVArProgr}
Let $\varepsilon > 0$ be an arbitrarily small constant, and $x \ge x_0(\varepsilon)$ a sufficiently large positive number. Let $q$ be a natural number such that $q \le x^{\kappa_0}$, where $\kappa_0$ is defined by~(\ref{formula-opred-ae0}). Suppose
\[
\alpha \in \left[\left(\theta + \varepsilon\right)\frac{\ln q}{\ln x}, 1 - 2.5\frac{\ln q}{\ln x}\right], \qquad x_1 \ge x^{1 - \alpha}, \qquad x_2 \ge x^\alpha,
\]
and $p_1$, $p_2$ are prime numbers. Assume $(a, q) = (l, q) = 1$, and let $\pi_2(x_1, x_2, a, l)$ denote the number of such pairs with $p_1 \le x_1$, $p_2 \le x_2$ and $p_1(p_2 + a) \equiv l \pmod q$. Then for any $A > 0$, the following asymptotic formula holds:
\[
\pi_2(x_1, x_2, a, l) = \frac{1}{\varphi(q)} \prod_{p \mid q} \left(1 - \frac{1}{p - 1}\right) \mathrm{Li}(x_1)\mathrm{Li}(x_2) + O\left(\frac{x_1 x_2}{\varphi(q) \ln x_1 \ln x_2 \mathscr{L}^A}\right),
\]
where the implied constant in $O$ depends only on $\varepsilon$.}

In the proof of Theorem~\ref{Teorema-RZKh-RPSdPrChVArProgr}, we also use the estimate~(\ref{formula-otsenka T(x;Q)-RZKh}), specifically its consequence (Lemma~\ref{Teorema-AAK-RPSdPrChVArProgr}) concerning the distribution of prime numbers in arithmetic progressions ``on average'', and the Brun--Titchmarsh theorem (Lemma~\ref{Lemma-Teorema Bruno-Titchmarsh}).

{\lemma\label{Lemma-Teorema Bruno-Titchmarsh}
For $(a, q) = 1$ and $q \le x$, we have
\[
\pi(x; q, a) \le \frac{2x}{\varphi(q)\ln\left(\frac{2x}{q}\right)}.
\]}

\section{Proof of Theorem~\ref{Teorema-RZKh-RPSdPrChVArProgr}}

\quad Without loss of generality, in view of the conditions of the theorem, we assume that the parameters $x_1$ and $x_2$ satisfy the following relations:
\begin{equation}\label{formula-usl-x1-i-x2}
q^\frac52\le x_1 \le q^{\frac52+c_1}, \qquad q^{\theta + \varepsilon} \le x_2 \le q^{\theta + \varepsilon + c_2},
\end{equation}
where $c_1$ and $c_2$ are arbitrary fixed positive constants. Hence,
\begin{equation}\label{formulas lnxi=lnq}
\ln x_1 \asymp \mathscr{L}, \qquad \ln x_2 \asymp \mathscr{L}.
\end{equation}
Using the orthogonality property of Dirichlet characters, we find:
\begin{align*}
\pi_2(x_1, x_2, a, l) &= \sum_{p_1 \le x_1} \sum_{p_2 \le x_2} \frac{1}{\varphi(q)} \sum_{\chi \bmod q} \chi(p_1(p_2 + a)) \overline{\chi}(l) \\
&= \frac{1}{\varphi(q)} \sum_{\chi \bmod q} T_1(x_2, \chi)\, \overline{\chi}(l) \sum_{p \le x_1} \chi(p).
\end{align*}
Splitting the last sum over $\chi$ into two parts, we obtain:
\begin{align}
\pi_2(x_1, x_2, a, l) &= \mathscr{M}_2(x_1, x_2, a, l) + \mathscr{R}_2(x_1, x_2, a, l), \label{formula H2=G2+R2+..} \\
\mathscr{M}_2(x_1, x_2, a, l) &= \frac{1}{\varphi(q)} T_1(x_2, \chi_0) \sum_{p \le x_1} \chi_0(p), \nonumber \\
\mathscr{R}_2(x_1, x_2, a, l) &= \frac{1}{\varphi(q)} \sum_{\chi \ne \chi_0} T_1(x_2, \chi)\, \overline{\chi}(l) \sum_{p \le x_1} \chi(p). \nonumber
\end{align}
In this formula, $\mathscr{M}_2(x_1, x_2, a, l)$ gives the expected main term of $\pi_2(x_1, x_2, a, l)$, while $\mathscr{R}_2(x_1, x_2, a, l)$ contributes to its error term.

Let us compute the main term. Using the asymptotic law of the distribution of prime numbers and relation~(\ref{formulas lnxi=lnq}), we obtain
\begin{equation}\label{formula-sump<=x1chi0(p)=x1/lnx1+O(...)}
\sum_{p \le x_1} \chi_0(p) = \pi(x_1) + O\left((\ln x_1)^2\right) = \mathrm{Li}(x_1) + O\left(\frac{x_1}{\ln x_1\, \mathscr{L}^A}\right).
\end{equation}
Now, represent the sum $T_1(x_2, \chi_0)$ in the form
\begin{align}
T_1(x_2, \chi_0) &= \sum_{\substack{p \le x_2 \\ (p+a,q)=1}} 1 = \mathrm{Li}(x_2) \sum_{d\mid q} \frac{\mu(d)}{\varphi(d)} + \sum_{d\mid q} \mu(d) \left( \pi(x_2; d, -a) - \frac{\mathrm{Li}(x_2)}{\varphi(d)} \right) \nonumber \\
&= \mathrm{Li}(x_2) \prod_{p\mid q}\left(1-\frac{1}{p-1}\right) + R_1(x_2) + R_2(x_2), \label{formula T(x2,chi=..+R1+R2} \\
R_1(x_2) &= \hspace{-15pt}\sum_{\substack{d\mid q \\ d \le \sqrt{x}(\ln x)^{-A-3.5}}} \hspace{-15pt}\left( \pi(x_2; d, -a) - \frac{\mathrm{Li}(x_2)}{\varphi(d)} \right), \nonumber \\
R_2(x_2) &= \hspace{-15pt}\sum_{\substack{d\mid q \\ d > \sqrt{x}(\ln x)^{-A-3.5}}} \hspace{-15pt}\left( \pi(x_2; d, -a) - \frac{\mathrm{Li}(x_2)}{\varphi(d)} \right). \nonumber
\end{align}
We now estimate $R_1(x_2)$ and $R_2(x_2)$ using, respectively, the Bombieri-Vinogradov theorem (Lemma~\ref{Sledst-teorema-RZKh-otsenka-T(x;Q)}) and the Brun-Titchmarsh theorem (Lemma~\ref{Lemma-Teorema Bruno-Titchmarsh}). We have:
\begin{align*}
R_1(x_2) &\ll \sum_{\substack{d\mid q \\ d \le \sqrt{x_2}(\ln x_2)^{-A-3.5}}} \left| \pi(x_2; d, -a) - \frac{\mathrm{Li}(x_2)}{\varphi(d)} \right| \ll \frac{x_2}{(\ln x_2)^{A+1}}; \\
R_2(x_2) &\le \sum_{\substack{d\mid q \\ d > \sqrt{x_2}(\ln x_2)^{-A-3.5}}} \hspace{-15pt}\left( \pi(x_2; d, -a) + \frac{\mathrm{Li}(x_2)}{\varphi(d)} \right) \ll \\
&\ll \sum_{\substack{d\mid q \\ d > \sqrt{x_2}(\ln x_2)^{-A-3.5}}} \hspace{-15pt}\left( \frac{x_2}{\varphi(d) \ln\left(\frac{x}{d}\right)} + \frac{x_2}{\varphi(d) \ln x_2} \right) \ll \\
&\ll \frac{x_2}{\ln x_2} \sum_{\substack{d\mid q \\ d > \sqrt{x_2}(\ln x_2)^{-A-2.5}}} \frac{1}{\varphi(d)} \ll \frac{x_2}{(\ln x_2)^{A+1}}.
\end{align*}
Substituting these estimates into~(\ref{formula T(x2,chi=..+R1+R2}) and using relation~(\ref{formulas lnxi=lnq}), we obtain
$$
T_1(x_2, \chi_0) = \prod_{p\mid q}\left(1-\frac{1}{p-1}\right) \mathrm{Li}(x_2) + O\left( \frac{x_2}{\ln x_2\, \mathscr{L}^A} \right).
$$
Thus, from~(\ref{formula-sump<=x1chi0(p)=x1/lnx1+O(...)}) and~(\ref{formula H2=G2+R2+..}), we find
\begin{align*}
\mathscr{M}_2(x_1, x_2, a, l) &= \frac{1}{\varphi(q)} \prod_{p\mid q} \left(1 - \frac{1}{p-1}\right) \mathrm{Li}(x_1) \mathrm{Li}(x_2) + O\left( \frac{x_1 x_2}{\varphi(q) \ln x_1 \ln x_2\, \mathscr{L}^A} \right).
\end{align*}

Let us estimate the remainder term $\mathscr{R}_2(x_1, x_2, a, l)$. Passing to estimates and then applying Lemmas~\ref{Teorema-RZKh-Otsenka-T(chi)-q-bezcub} and~\ref{Teorema-RZKh-Otsenka-T(chi)-q-sost} to estimate the sum $T_1(x_2,\chi)$ for $\chi \neq \chi_0$, and using Lemma~\ref{Sled-teorema-RZKh-otsenka-t(x;q)}, we obtain
\begin{align*}
\mathscr{R}_2(x_1, x_2, a, l) &\le \frac{1}{\varphi(q)} \max_{\chi\neq\chi_0} |T_1(x_2,\chi)| \sum_{\chi \bmod q} \left| \sum_{p\le x_1} \chi(p) \right| \ll \\
&\ll \frac{1}{\varphi(q)} x_2 \exp\left( -0.6 \sqrt{\mathscr{L}} \right) \left( x_1 (\ln x_1 q)^2 + x_1^{\frac{4}{5}} q^{\frac{1}{2}} (\ln x_1 q)^{33} + x_1^{\frac{1}{2}} q (\ln x_1 q)^{33} \right) = \\
&= \frac{1}{\varphi(q)} \cdot \frac{x_1 x_2}{\ln x_1 \ln x_2} \Delta(x_1, x_2, q),
\end{align*}
where
\begin{align*}
\Delta(x_1, x_2, q) &= \left( (\ln x_1 q)^{-31} + x_1^{-\frac{1}{5}} q^{\frac{1}{2}} + x_1^{-\frac{1}{2}} q \right) (\ln x_1 q)^{33} \ln x_1 \ln x_2 \exp\left( -0.6 \sqrt{\mathscr{L}} \right).
\end{align*}

Thus, from relations~(\ref{formula-usl-x1-i-x2}) and~(\ref{formulas lnxi=lnq}), we have
\begin{align*}
\Delta(x_1, x_2, q) &\ll (\ln x_1 q)^{33} \ln x_1 \ln x_2 \exp\left( -0.6 \sqrt{\mathscr{L}} \right) \ll \mathscr{L}^{35} \exp\left( -0.6 \sqrt{\mathscr{L}} \right) \ll \mathscr{L}^{-A}.
\end{align*}

The theorem is thus proved.



\vspace{10mm}
\begin{thebibliography}{12} 
\bibitem{Karatsuba-DANSSSR-1970-192-4}  Karatsuba~A.~A., 1970, ``The distribution of products of shifted prime numbers in arithmetic progressions'', \emph{Dokl. Akad. Nauk SSSR}, vol.~192, Is.~4, pp.~724-727.
\bibitem{KaratsubaAA-2008} Karatsuba,~A.~A., 2008, ``Arithmetic problems in the theory of Dirichlet characters'', \emph{Russian Mathematical Surveys}, vol~63, Is.~4, pp.~641-690.
\bibitem{PetechukMM-IzvFNSSSR-1979-43-4}  Petecuk~M.~M., 1980,``The sum of the values of the divisor function in arithmetic progressions whose difference is a power of an odd prime'', \emph{Mathematics of the USSR-Izvestiya}, vol.~15, Is.~1, pp.~145-160.
\bibitem{Chubarikov-Vestnik-MGU-1973-2}  Chubarikov~V.~N., 1973, ``A more precise bound for the zeros of {Dirichlet} $L$-series modulo with a power of prime'', \emph{Moscow University Mathematics Bulletin}, vol.~28, no~1-2, pp.~76-81.
\bibitem{Friendlander+Iwaniec-ActaArith-1985-45-3}  Friendlander~J.~B., \& Iwaniec~H., 1985, ``The divisor problem for arithemetic progressions'', \emph{Acta Arith.}, vol.~45, Is.~3, pp.~273-277.
\bibitem{RakhmonovZKh-Izv-ANSSSR-1989}  Rakhmonov~Z.~Kh., 1990, ``The distribution of Hardy-Littlewood numbers in arithmetic progressions'', \emph{Mathematics of the USSR-Izvestiya}, vol.~34, Is.~1, pp.~213-228. %
\bibitem{PanChD+HfnChB-OATCh-1991-Pekin} Pan Chengdong, \& Pan Chengbiao, 1991,``Foundation to Analytic Number Theory'', \emph{Science Press, Beijing, 1991}, (in Chinese).
\bibitem {Montgomeri-1974} Montgomery,~H., 1971, \emph{Topics in Multiplicative Number Theory}, vol. 227. Springer-Verlag, Berlin-New York.
\bibitem {Vaughan-1975} Vaughan,~R.~O., 1975, ``Mean value theorems in prime number theory'', \emph{J.London Math. Soc.}, vol.~s2-10, Is.~2, pp.~153-162.
\bibitem{RakhmonovZKh-IzvRAN-1993-57-4} Rakhmonov,~Z.~Kh., 1994,~``Theorem on the mean value of $\psi(x,\chi)$ and its applications'', \emph{Russian Academy of Sciences. Izvestiya Mathematics}, vol.~43, Is.~1, pp.~49-64.
\bibitem{RakhmonovZKh-DAN-1993-331-3} Rakhmonov,~Z.~Kh., 1994,~``Mean values of the Chebyshev function'', \emph{Russ. Acad. Sci., Dokl., Math.}, vol.~48, Is.~1, pp.~85-87.
\bibitem{RZKh+NOO-ChebSb-2021-22-59(81)}  Rakhmonov,~Z.~Kh., \& Nozirov~O.~O., 2021,~  ``On the mean values of the Chebyshev function and their applications'',  \emph{Chebyshevskii Sbornik}, vol.~22, no~5(81), pp.~198-222.
\bibitem{RakhmonovZKh-1994-IzvRAN} Rakhmonov,~Z.~Kh., 1995,~``A mean-value theorem for Chebyshev functions'', \emph{Russian Academy of Sciences. Izvestiya Mathematics}, vol.~44, Is.~3, pp.~555-569.
\bibitem{RakhmonovZKh-DAN-1996} Rakhmonov,~Z.~Kh., 1996,~``The mean-value theorem in prime number theory'', \emph{Doklady Mathematics}, vol.~54, Is.~1, pp.~597-598.
\bibitem{Timofeev-Izv-ANSSSR-1987} Timofeev,~N.~M., 1988,  ``Distribution in the mean of arithmetic functions in short intervals in progressions'', \emph{Mathematics of the USSR-Izvestiya}, vol.~30, Is.~2, pp.~315-335.
\bibitem{VinigradovIM-1938} Vinogradov, I.~M., 1938, ``On the distribution of quadratic rests and non-rests of the form $p+k$ to a prime modulus'', \emph{Rec. Math. [Mat. Sbornik] N.S.}, vol.~3(45), no~2,  pp.~311-319.
\bibitem{VinigradovIM-1943} Vinogradov,~I.~M., 1943, ``An improvement of the estimation of sums with primes'', \emph{Izv. Akad. Nauk SSSR. Ser. Mat.}, vol.~7, no~1, pp.~17-34.
\bibitem{Jutila} Jutila,~M., 1968, ``On the least Goldbach's number in an arithmetical progression with a prime difference'', \emph{Ann. Univ. Turku; Ser. A,} I 118(5).
\bibitem{VinigradovIM-1952}  Vinogradov,~I.~M., 1952, ``New approach to the estimation of a sum of values of $\chi(p+k)$'', \emph{
Izv.~Akad.~Nauk~SSSR. Ser.~Mat.}, vol.~16, Is.~3, pp.~197-210.
\bibitem{VinigradovIM-1953}  Vinogradov,~I.~M., 1953, ``Improvement of an estimate for the sum of the values  $\chi(p+k)$'', \emph{Izv. Akad. Nauk SSSR. Ser. Mat.}, vol.~17, Is.~4, pp.~285-290.
\bibitem{Linnik-1973}  Linnik,~Yu.~V., 1975,  ``Recent works of I.M. Vinogradov'', \emph{Proceedings of the Steklov Institute of Mathematics},  vol.~132, pp.~25-28.
\bibitem{KaratsubaAA-1968} Karatsuba,~A.~A., 1968, ``Sums of characters, and primitive roots, in finite fields'', \emph{Dokl. Akad. Nauk SSSR}, vol.~180, Is.~6.  pp.~1287-1289.
\bibitem{KaratsubaAA-1970-1} Karatsuba,~A.~A., 1970, ``Estimates of character sums'', \emph{Math. USSR-Izv.}, vol.~4, Is.~1, pp.~19-29.
 \bibitem{KaratsubaAA-1970-2} Karatsuba,~A.~A., 1970, ``Sums of characters over prime number'', \emph{Math. USSR-Izv.}, vol.~4, Is.~2, pp.~303-326.
\bibitem{KaratsubaAA-1970-3} Karatsuba,~A.~A., 1970, ``Sums of characters with prime numbers'', \emph{Dokl. Akad. Nauk SSSR}, vol.~190, Is.~3, pp.~517-518.
\bibitem{RakhmonovZKh-1986-UMN} Rakhmonov,~Z.~Kh.,~1986,~``On the distribution of values of Dirichlet characters'', \emph{Russian Math. Surveys}, vol.~41, Is.~1, pp/~237--238.
\bibitem{RakhmonovZKh-1986-DANRT} Rakhmonov,~Z.~Kh.,~1986, ``Estimation of the sum of characters with primes'', \emph{Dokl. Akad. Nauk Tadzhik. SSR}, vol.~29, Is.~1, pp.~16--20,, (in Russian).
\bibitem{RakhmonovZKh-1994-TrMIRAN} Rakhmonov,~Z.~Kh., 1995,~``On the distribution of the values of Dirichlet characters and their applications'', \emph{Proc. Steklov Inst. Math.},  vol.~207,  pp.~263--272.
\bibitem{RakhmonovZKh-1986-IzvANRT} Rakhmonov,~Z.~Kh., 1986~``The least Goldbach number in an arithmetic progression'', \emph{Izv. Akad. Nauk Tadzhik. SSR. Otdel. Fiz.-Mat., Khim. i Geol. Nauk}, Is.~2(100), pp.~103-106, (in Russian).
\bibitem{Huxley-1942-Inv.Math} Huxley,~M.~N., 1971, ``On the difference between consecutive primes'', \emph{Inventiones mathematicae}, vol.~15, Is.~2, pp.~164-170.
\bibitem{Friedlander-Gong-Shparlinski}Fridlander,~Dzh.~B., \& Gong,~K., \& Shparlinskii,~I.~E., 2010, ``Character sums over shifted primes'', \emph{Math. Notes}, vol.~88, Is.~3-4, pp.~585-598.
\bibitem{RakhmonovZKh-2013-DANRT} Rakhmonov,~Z.~Kh., 2013,~``Distribution of values of Dirichlet characters in the sequence of shifted primes'', \emph{Doklady Akademii nauk Respubliki Tajikistan},  vol.~56, Is.~1, pp.~5-9, (in Russian).
\bibitem{RakhmonovZKh-2013-IzvSarUniv} Rakhmonov,~Z.~Kh., 2013,~``Distribution of values of Dirichlet characters in the sequence of shifted primes'', \emph{Izv. Saratov Univ. (N.S.), Ser. Math. Mech. Inform.}, vol.~13, Is.~4(2), pp.~113-117, (in Russian).
\bibitem{RakhmonovZKh-2014-Ch.sbor-15-2(50)} Rakhmonov,~Z.~Kh.,~2014, ``Sums of characters over prime numbers'', \emph{Chebyshevskii Sb.}, vol.~15, no~2, pp.~73-100, (in Russian).
\bibitem{Kerr-2021} Kerr,~B., 2021,~``Bounds of Multiplicative Character Sums over Shifted Primes'',  \emph{Proc. Steklov Inst. Math.}, vol.~314, pp.~64-89. 
\bibitem{RakhmonovZKh-2017-TrMIRAN-299} Rakhmonov,~Z.~Kh.,~2017, ``Sums of values of nonprincipal characters over a sequence of shifted primes'', \emph{Proc. Steklov Inst. Math.}, vol.~299, pp.~219-245.
\bibitem{RakhmonovZKh-2017-DANRT-9} Rakhmonov,~Z.~Kh., 2017,~``On the estimation of the sum the values of Dirichlet character in a sequence of shifted primes'', \emph{Doklady Akademii nauk Respubliki Tajikistan}, vol.~60, no~9, pp.~378-382, (in Russian).
\bibitem{RakhmonovZKh-2025-SummNeglChar-arxiv}	Rakhmonov~Z., 2025, ``Sums   of   values   of   non-principal  characters  over shifted primes'', \emph{https://arxiv.org/abs/2412.03040}, 12 March, 2025,  doi.org/10.48550/arXiv.2412.03040.
\bibitem{RakhmonovZKh-2018-Springer} Rakhmonov~Z.~Kh., 2018, ``Sums of Values of Nonprincipal Characters over Shifted Primes'',  In: Pintz J., Rassias M. (eds) Irregularities in the Distribution of Prime Numbers, \emph{Springer International Publishing}, pp.~187-217.
\end{thebibliography}
\end{document}